\documentclass[12pt]{amsart}
\theoremstyle{plain}

 \newtheorem{Thm}{Theorem}[section]
 
 \newtheorem{Lem}[Thm]{Lemma}
 \newtheorem{Prop}[Thm]{Proposition}
 \newtheorem{Def}[Thm]{Definition}
 \theoremstyle{Remark}
 
 \numberwithin{equation}{section}

\begin{document}

\title[subelliptic equations]
{Some results for subelliptic equations}

\author{Li Ma$^*$, Dezhong Chen, and Yang Yang}

\address{Ma: Department of Mathematical Sciences\\
Tsinghua University\\ Beijing 100084, P. R. CHINA}
\email{lma@math.tsinghua.edu.cn}
\address{ Chen and Yang: Department of Mathematical Sciences\\
Tsinghua University\\ Beijing 100084, P. R. CHINA}
\date{}

\subjclass{Primary 35J65, 35H20}

\keywords{principal eigenvalue, Hormander's Laplacian, positive
solution}

\thanks{$^*$ This work is supported in part by
the key 973 project of China. Part of this work was done while the
first named author was working in NUS in Singapore in the year 2001.
The first named author would also like to thank Prof.Hua Chen for
inviting him to give lectures about sub-elliptic problems in Wuhan
University}

\begin{abstract}
In this paper, we consider the principal eigenvalue problem for
Hormander's Laplacian on $R^n$, and we find a comparison principle
for such principal eigenvalues. We also study a related semi-linear
sub-elliptic equation in the whole $R^n$ and prove that under a
suitable condition, we have infinite many positive solutions of the
problem.
\end{abstract}

\maketitle

\section{Introduction}

In this work, we study sub-elliptic problems arising from
Hormander's vector fields on the whole Euclidean space $R^n$. We
shall study the principal eigenvalue problem of Hormander's
Laplacian on $R^n$.  We shall also treat some related semi-linear
sub-elliptic problems both in bounded domains of $R^n$ and in $R^n$.
Recently, people studied the geometry and analysis for the
Hormander's vector fields, and some basic regularity properties of
Hormander's Laplacian on bounded or closed manifolds have been found
(see \cite{CGL}, \cite{Gr}, \cite{X}). In fact, many important
geometric problems (such as the Yamabe problem on CR manifolds and
sub-elliptic harmonic maps) have attracted a lot of attentions (see
\cite{JL} and \cite{JX}). From the analytical point of view, it is a
challenge to extend some of the beautiful results for elliptic
problems to those of sub-elliptic problems. Interestingly, the
moving plane method has been used to study the symmetry properties
of positive entire solutions of Yamabe-type equations on groups of
Heisenberg type (see \cite{BP}, \cite{GV}, \cite{LW}). All these
extensions are nontrivial. The potential mathematical applications
of this research direction may be found in complex geometry.

Before giving the statements of our results, let us now recall some
definitions and properties concerning with Hormander's vector fields
on $R^n$.

Given a family of smooth vector fields $X_j$, $1\leq j\leq m$ in
$R^n$. We say that $\{X_j\}$ are {\bf Hormander's vector fields} if
$X_1,...,X_m$ together with their communicators up to a certain
fixed length $p$ span the tangent space $R^n$ at each point $x\in
R^n$. The important thing here is that from this system of vector
fields, we can define a metric in the following way. Let
$\mathbb{H}$ be the space of horizontal curves, i.e., each such
curve is a piecewise $C^1$ mapping
$$
\gamma:[0,T]\to R^n
$$
such that it holds
$$
\gamma'(t)=\sum_{j=1}^m f_j(t)X_j(\gamma(t)).
$$
for some piecewise continuous functions $f_1(t),...,f_m(t)$ with
$$
\sum_{j=1}^m f_j^2(t)\leq 1.
$$
We call $T$ the horizontal length of $\gamma$. For any $x,y\in
R^n$, define
$$
d(x,y)=inf_{\gamma\in\mathbb{H}}\{T; \gamma(0)=x,\gamma(T)=y\}
$$
It is known that this metric is equivalent to the usual Euclidean
metric (\cite{NSW}\cite{HK}). We denote by $B_d(x_0,R)$ the metric
ball in $(R^n,d)$ with center $x_0\in R^n$ and radius $R>0$.

Let $D$ be an open domain in $R^n$. We introduce the following
Sobolev spaces. For $f\in C^1_0(D)$, let
$$
|Xf(x)|^2=\sum_j|X_jf(x)|^2,
$$
where $x\in D$. We define
$$
|Xf|_{M^1}=(\int_{D}|Xf(x)|^2dx)^{1/2},
$$
which is a norm defined on $C^1_0(D)$, and let $M^1_0(D)$ be the
completion of $C^1_0(D)$ under the norm $|\cdot|_{M^1}$. We let
$$
M^1(D)=\{f\in L^2(D);|Xf|\in L^2(D)\}.
$$

We also use the definition
$$
M^1_{loc}=\{f\in L^1_{loc}(R^n); \int_{R^n}|Xf(x)|^2dx<+\infty \}.
$$
for the whole space $D=R^n$.

  Write by
  $$
X_j(x)=A^{jk}(x)\partial_{x^k}.
  $$
  Here and hereafter we use the standard summation convention.
  Throughout this paper,
   we always assume that
 the vector field group  $\{X_j\}_1^m$ on $R^n$ is
 a family of Hormander's vector fields.
For each $j$, we define the formal adjoint of $X_j$ by
$$
X^*_jf=-\partial_{x^k}(A^{jk}(x)f),
$$
where $f\in C^1_0(R^n)$. Then we give the definition of
\emph{Hormander's Laplacian }on $R^n$ (see also \cite{X}):
$$
H=\sum_j
X^*_jX_j=-\sum_{i,k}\partial_{x^i}(a^{ik}(x)\partial_{x^k}),
$$
where
$$
a^{ik}(x)=\sum_jA^{ji}(x)A^{jk}(x).
$$
{\bf Remarks}. Generally speaking, this operator is not elliptic,
but it enjoys a lot of properties from elliptic operators. In
particular, $H$ is a hypoelliptic operator (see \cite{H}). We also
have many properties just like in the elliptic case. For example,
$M^1_0(D)$ is a Hilbert space for bounded domain $D\subset R^n$, and
we have Bony's maximum principle (see \cite{B}), the Harnack
inequality (see \cite{CGL}), and Sobolev's inequality (see
\cite{GN}) under suitable conditions. Therefore, people can use
variational methods to prove the existence of weak solutions of some
semi-linear partial differential equations
 related to Hormander's Laplacian. Some well-known results will be recalled in
 the next section. In this paper, we say that a bounded domain $D$ is
 a {\em regular bounded domain} if its boundary $\partial D$ has
 non-characteristic points for the Hormander's Laplacian $H$.

Given a smooth function $V(x)$, which may change sign in $R^n$.
People want to study the unique continuation property of the
following equation
\begin{equation}
Hu-Vu= 0, \quad{ in}\quad R^n. \label{HV}
\end{equation}

We will use another type of unique continuation property \cite{Ga}.

\begin{Def} \label{1.1} We say that {\em the weak unique continuation
property} is true for equation (\ref{HV}) if for every bounded
regular domain $\Omega$, $u\in M^1_0(\Omega)$ is a non-negative weak
solution of
$$Hu-Vu= 0, \quad{ in}\quad \Omega $$
and $u=0$ on an open subset of $\Omega$, then $u=0$ in $\Omega$.
\end{Def}

Harnack Inequality for (\ref{HV}) has a closed relation with unique
continuation property. Hence, we assume that there holds {\em the
Harnack inequality} for non-negative solutions of the equation
(\ref{HV}) on $\Omega$. Namely, if $u$ is a non-negative smooth
solution of the equation (\ref{HV}) on $\Omega$, then, for each
point $x_0\in \Omega$, there exist positive constants $R$ and $C$
such that
$$
\sup_{B_d(x_0,R)}u\leq C\inf_{B_d(x_0,R)}u.
$$

{\bf Remarks}. We point out that our {\em weak unique continuation
property} follows from the Harnack inequality. Note that this
Harnack inequality is true for smooth potentials $V$ (see
\cite{CGL}). In fact, a larger class of potentials $V$ for the
Harnack inequality being true is found there. Assume that the vector
field group $X_j$, $1\leq j\leq m$ is a family of Hormander's vector
fields on $R^n$. In \cite{CGL}, the authors found the class
$\mathcal{K}_H^{loc}$ for the potential functions $V$ such that the
Harnack inequality is true for the non-negative solutions of
(\ref{HV}). Let's recall the definition for the interested readers.
Let $\Gamma$ be a positive fundamental solution of (\ref{HV}) which
is $C^{\infty}$ off the diagonal in $R^n\times R^n$. For $x\in R^n$,
$r>0$, let
$$
\Omega_r(x)=\{y\in R^n;\Gamma(x,y)>r^{-1}\}.
$$
Then $V\in \mathcal{K}_H^{loc}$ if
$$
\lim_{r\to 0} \sup_{x\in U} \int_{\Omega_r(x)}|V(y)|\Gamma(x,y)dy
=0 $$ for every bounded set $U\in R^n$.

With these preparations, we state our following results.

\begin{Thm}\label{T1.2}  Given a smooth function $V$ on $R^n$. If there is a
positive smooth function $u$ satisfying $Hu-Vu\geq 0$ in $R^n$, then
for any bounded regular domain $D\subset R^n$, we have
$$
\lambda^H_1(D):=inf\{\int_D(|Xf|^2-Vf^2)dx; f\in C^1_0(D),
\int_Df^2dx=1\}> 0.
$$
\end{Thm}

We remark that, although this is a new result for sub-elliptic
operator, its elliptic version is well-known in the theory of
elliptic partial differential equations of second order (see
\cite{FS}). Actually, the proof of Theorem \ref{T1.2} is almost the
same as in elliptic problems. However, for completeness, we will
give a detailed proof for our case in the next section. In our proof
of Theorem \ref{T1.2}, we only make use the interior regularity of
weak solution and our weak unique continuation property. It would be
interesting to study viscosity solutions of non-linear heat
flow/system with Hormander's Laplacian (see \cite{LL} for related
work). We also point out that every $\lambda^H_1(D)$ is achieved in
the sense that there is a positive smooth function $u_D$ satisfying
$$ Hu-Vu= \lambda^H_1(D)u $$
  with $u_D=0$ on $\partial D$. One can use the direct method and the compact
  imbedding theorem (see Theorem \ref{T2.2} below) to obtain this result.
  Another way to get this result is the approximation method. We give an idea
  for approximation method here. Define
  $$
  J(f)=\int_D(|Xf|^2-Vf^2)dx=\int_D(a^{ik}\partial_if\partial_kf-Vf^2)dx
  $$
for $f\in C^1_0(D)$. By Lemma 7.6 in \cite{GT}, we have
 $\frac{\partial}{\partial x^i} |f|= \frac{\partial}{\partial x^i}
  f$ in $\{x;f(x)>0\}$, $\frac{\partial}{\partial x^i} |f|=0$ in $\{x;f(x)=0\}$,
and $\frac{\partial}{\partial x^i} |f|= -\frac{\partial}{\partial
x^i} f$ in $\{x;f(x)<0\}$. Then
  $J(|f|)=J(f)$ for $f\in C^1_0(D)$. For $\epsilon>0$, we let
$$
J_{\epsilon}(f)=\int_D((a^{ik}+\epsilon\delta_{ik})\partial_if\partial_kf-Vf^2)dx
$$
Then we also have that $J_{\epsilon}(|f|)=J_{\epsilon}(f)$ for
$f\in C^1_0(D)$. Define
$$
\lambda^{\epsilon}_1(D):=inf\{J_{\epsilon}(f); f\in C^1_0(D),
\int_Df^2dx=1\}> 0.
$$
Then $\lambda^{\epsilon}_1(D)$ is strictly decreasing to
$\lambda^H_1(D)$ as $\epsilon\to 0+$.

 Hence, by standard elliptic theory \cite{GT}, we can find a positive smooth function
 $u_{\epsilon}\in C^2_0(D)$ satisfying
 $$
(H-\epsilon
\Delta-V)u_{\epsilon}=\lambda^{\epsilon}_1(D)u_{\epsilon}
 $$
with $J_{\epsilon}(u_{\epsilon})=\lambda^{\epsilon}_1(D)$.

Using a priori estimate (see Lemma 3.6 in \cite{JX}) we can find a
subsequence
  $(u_k)\subset (u_{\epsilon})$ which converges in $M^2_{loc}(D)$ to a nontrivial
  non-negative function $u$. This $u$ satisfies
 $$ Hu-Vu= \lambda^H_1(D)u $$
 and
$$
\lim_{k\to +\infty}J(u_k)=\lambda^H_1(D).
$$
Using the regularity result ( see \cite{HK} and \cite{JL}) we know
that $u$ is the eigenfunction of $\lambda^H_1(D)$.

We come to consider the principal eigenvalue comparison results. Let
$g(x)$ be a nontrivial (maybe non-smooth) function, which may change
sign in $R^n$.
 We consider the principal eigenvalue problem for the
equation:
\begin{equation}
Hu=\lambda g(x)u, \quad{in}\quad{R^n.} \label{g}
\end{equation}
By definition, {\em a principal eigenvalue} of (\ref{g}) is a
positive constant $\lambda_0$ such that there is a positive $C^1$
solution $u(x)$ of (\ref{g}) when $\lambda=\lambda_0$. When
$H=-\Delta$ is the standard positive Laplace operator on $R^n$,
there have been a number of authors studying such a problem. One may
see \cite{BNV}  for references. However, there is few work on our
problem (\ref{g}). We have the following result

\begin{Thm} \label{T1.3} Assume that $g_+ $ is a smooth potential (or more
generally $g_+\in \mathcal{K}_H^{loc}$) such that there is a
principal eigenvalue $\mu$ for the problem (\ref{g}) with $g=g_+$.
Let $g(x)$ be a smooth function in $R^n$ such that
$$
g(x)\leq g_+(x), x\in R^n.
$$
Then any number $\lambda\in (0, \mu]$ is a principal eigenvalue of
(\ref{g}).
\end{Thm}

{\bf Remarks}. We believe that this kind of result is also true for
some non-linear sub-elliptic problems. In the standard elliptic
case, this result has been obtained by Z.R. Jin in his interesting
paper \cite{J}. Clearly, our result is more general than his result.
Jin's argument can not be carried out to our case because he used
stronger estimates for second order uniformly elliptic operators.
Our argument is new and it is of variational nature. Our proof is
based on Theorem \ref{T1.2}.

We also study the existence of a non-trivial non-negative solution
of the following sub-elliptic equations on regular domains:
$$
   -H u+\lambda a(x)u-b(x)u^p=0 \quad{in}\quad \Omega
$$
with the Dirichlet boundary condition:
$$
u=0\quad{on}\quad \partial\Omega,
$$
where $a(x)$ and $b(x)$ are continuous functions on $\Omega$ with
$\|a\|_\infty<\infty$ and $b(x)$ nonnegative and not identically
zero, $\Omega$ is a smooth bounded regular domain in $R^n$. our
result is stated in Proposition \ref{4.2} below.

As we said before, we are also interested in studying the existence
of positive solutions of a semi-linear sub-elliptic equation in the
the whole space $R^n$:
\begin{equation}
Hu+k(x)u-K(x)|u|^{p-1}u=0 \quad{in}\quad R^n,\label{Yam}
\end{equation}
where $k(x)$ and $K(x)$ are given smooth functions in $R^n$, $p>1$
is a constant. It is clear that the equation is a Yamabe-type or
logistic equation. To get a reasonable result, we need a basic
property for the vector fields $\{ X_j\}$ in $R^n$.

{\bf Property (P)}: \emph{ There exists a suitable non-negative
function $f\in C^{\alpha}_{loc}(R^n)$ such that for given any
$\epsilon \in {R}$, we have the property that the Poisson equation
$$
HU=f, \quad{in}\quad R^n
$$
has a unique bounded solution $U\in C^{\alpha}(R^n)\cup
M^1_{loc}(R^n)$ such that $$\lim_{|x|\to +\infty}U(x)=\epsilon.$$}

{\bf Remark}. Note that Property $(P)$ is always true when $H$ is
uniformly elliptic in $R^n$. Property $(P)$ is also true for the
Hormander-Laplacian on the Heisenberg group (\cite{BCC},
\cite{JL}). We do not know if Property (P) is true for a general
class of Hormander vector fields.

Then we have the following

\begin{Thm} \label{T1.4}  Assume that the property (P) is true for
the Hormander vector fields $X_j$ in $R^n$.  Then there is a
constant $\theta>0$ such that if $|k(x)|\leq \theta f(x)$ and
$|K(x)|\leq \theta f(x)$ for every $x\in R^n$, then (\ref{Yam}) has
a family positive solutions $U_{\epsilon}$ in $ C^{\alpha}(R^n)\cup
M^1_{loc}(R^n)$ with $\lim_{|x|\to
+\infty}U_{\epsilon}(x)=\epsilon$.
\end{Thm}

We remark that a similar result was obtained by F. H. Lin (\cite{L})
when $H$ is uniformly elliptic in $R^n$. We think that such result
may play an important role in understanding sub-elliptic Yamabe-type
problems and Logistic equations on $R^n$.

Here is the plan of the paper. We shall state some well known facts
in section two. Proofs of Theorems \ref{T1.2} and \ref{T1.3} are
given in section three. We study the non-linear sub-elliptic problem
with the Dirichlet boundary condition in bounded domains in section
four. In section five, we study the Yamabe type or the logistic type
nonlinear problem (\ref{Yam}).

In the following, we denote $C$ the varying constants in different
positions.

\section{Some facts about Hormander's vector fields}

We first state Bony's maximum principle \cite{B}. Let $\Omega $ be a
regular domain in $R^n$ and let $\rho(x,y)$ be the
  distance function defined by the Hormander vector fields
  $\{X_j\}$. For $\alpha \in (0,1)$, we define the Holder spaces
  $$
S^{\alpha}(\Omega)=\{u\in L^{\infty}(\Omega);
[u]_X^{\alpha}=sup_{x,y\in \Omega}
\frac{|u(x)-u(y)|}{\rho(x,y)^{\alpha}}<+\infty\}.
  $$
  and
  $$
S^{1,\alpha}(\Omega)=\{u\in S^{\alpha}(\Omega); Xu\in
S^1(\Omega)\}.
  $$

\begin{Prop} \label{T2.1}(Bony's maximum principle): Assume that $w\in
S^{1,\alpha}(\Omega)$ satisfies
$$
(H-V)w\leq 0,
$$
in weak sense. If there is a point $x_0\in \Omega$ with
$$
0=w(x_0)=\sup_{\Omega} w(x),$$ then $w=0$ in $\Omega$. Furthermore,
the maximum principle also holds for weak solutions $w\in
M^1(\Omega)$ of $(H-V)w\leq 0$ in the usual weak sense .
\end{Prop}

A remarkable result for Hormander's vector fields is the following
 Poincare inequality of Jerison (see Theorem 11.20 and the remark afterward
 in \cite{HK}):

 \begin{Thm} \label{T2.2} Let $X_1,...,X_m$ be a system of Hormander's vector fields
 on $R^n$. Then for every compact subset $K\subset R^n$ there are
 constants $C$ and $R_0$ such that for any Lipschitz function $u:B_d\to R$
$$
\int_{B_d}|u-u_B|dx\leq CR\int_{B_d}|Xu|dx,
$$
where $B_d$ is a ball center at $K$ with radius $R<R_0$ and $u_B$
is the average of the integral of $u$ over $B_d$.
 \end{Thm}

The doubling condition (in short,DC) is true for Hormander's vector
fields on $R^n$ (\cite{NSW}): For any bounded open domain $D$ of
$R^n$, there exist positive constants $C_1$ and $R_1$ such that for
any $x_1\in D$, $R\leq R_1$, such that
$$
|B_d(x_1,2R)|\leq C_1 |B_d(x_1,R)|.
$$
Here $|\cdot|$ denotes the Lebesgue measure of the subset $\cdot$
of $R^n$.

The property (DC) is very important. With this property, one can
show that the Poincare inequality (see page 79 in \cite{HK}) implies
the following Sobolev inequality:
$$
(\frac{1}{|B_d|}\int_{B_d}|u|^qdx)^{1/q}\leq
CR(\frac{1}{|B_d|}\int_{B_{d}}|Xu|^pdx)^{1/p},
$$
with $q>p$ for all smooth functions $u$ with compact support in
the metric ball $B_d$.

By using these Poincare and Sobolev inequalities, we have the
Rellich-Kondrachov compactness theorem (see Theorem 8.1 in
\cite{HK}):
\begin{Thm}\label{T2.3}
Let $D$ be a bounded Lipschitz domain in $R^n$. Then the imbedding
from $M^1(D)$ into $L^2(D)$ is compact.
\end{Thm}

\section{Proofs of Theorems \ref{T1.2} and \ref{T1.3}}

 We now prove  Theorem \ref{T1.2}.
 \begin{proof}
Assume that we have a positive smooth function $u$ satisfying
$$
Hu-Vu \geq 0, \quad{in}\quad R^n
$$
 Let $w=log u$. We calculate and find:
\begin{equation}
Hw-|Xw|^2-V\geq  0. \label{3.1}
\end{equation}
Let $f\in C^1_0(D)$. Multiplying (\ref{3.1}) by $f^2$ and
integrating by parts, we have
$$
\int |Xw|^2f^2+Vf^2\leq \int f\langle Xf,Xw\rangle.
$$
By the Cauchy-Schwartz inequality we know that the right hand side
of the above inequality is bounded by
$$
\int f^2|Xw|^2+|Xf|^2.
$$
Then we get
$$
\int |Xf|^2-Vf^2\geq 0.
$$
This gives that $\lambda^H_1(D)\geq 0$. Using the weak unique
continuation property, this $\lambda^H_1(D)\geq 0$ is clearly
equivalent to the fact
$$
\lambda^H_1(D)> 0.
$$
In fact, if $\lambda^H_1(D)=0$ for some bounded regular domain $D$,
then $\lambda^H_1(\hat{D})=0$ for every larger regular domain
$\hat{D}$ containing $D$. Let $u_D$ be the principal eigenfunction
for $\lambda^H_1(D)=0$ and we extend $u_D$ to the domain $\hat{D}$
by zero extension. Then $u_D$ is also the eigenfunction for
$\lambda^H_1(\hat{D})=0$. Using the regularity theory of Hormander's
Laplacian operator (see \cite{JX} and \cite{X}), we know that $u_D$
is smooth in the interior of the domain $\hat{D}$. By our weak
unique continuation property we find that $u_D=0$ in $\hat{D}$. This
is a contradiction to the fact that $u_D>0$ in $D$. Thus, we have
proved Theorem \ref{T1.2}.
\end{proof}
\bigskip

In the rest of this section, we give the proof of Theorem
\ref{T1.3}.

\begin{proof} Take a fixed $\lambda\in (0, \mu]$. We
choose a sequence of bounded domain ${D_k}$ such that $0\in
D_k\subset D_{k+1}$ and $R^n=\bigcup_k D_k$ . We use the direct
method \cite{Ma} to solve the following problem
\begin{equation}
Hu= \lambda g(x)u, \quad{in}\quad D_k \label{3.2}
\end{equation}
with the boundary condition
$$
u=k,\quad {on}\quad \partial D_k.
$$

Fix $D=D_k$. Define the functional
$$
J(u)=\int_D |Xu|^2-\lambda g(x)u^2
$$
on the convex closed subset $E:=\{u=v+k; v\in M^1_0(D)\}$. Note
that $J$ is bounded from below on $E$. In fact, Let $I(v)=J(v+k)$.
Then
$$
I(v)=\int_D |Xv|^2-\lambda g(x)v^2-2k\lambda g(x)v-k^2\lambda g(x)
$$
Since
$$
\lambda g(x)v^2\leq \mu g_+(x)v^2,
$$
we have
$$
I(v)\geq\int_D |Xv|^2-\mu g_+(x)v^2-2k\lambda g(x)v-k^2\lambda
g(x)
$$
By the assumption of Theorem \ref{T1.3} and using Theorem
\ref{T1.2}, we have
$$
\int_D |Xv|^2-\mu g_+(x)v^2\geq \lambda^H_1(D)\int_Dv^2.
$$
So we have
$$
I(v)\geq \lambda^H_1(D)\int_Dv^2-\int_D 2k\lambda g(x)v-k^2\lambda
g(x)
$$
By the Cauchy-Schwartz inequality we have
$$ |2k \lambda\int_D
g(x)v|\leq \frac{1}{2}\lambda^H_1(D)\int_Dv^2+C\int_Dg^2
$$
Hence, we have
$$
inf_EJ(u)\geq -k^2\int \lambda g(x)-C\int_Dg^2.
$$
As we showed in the remarks after Theorem \ref{T1.2}, we can take a
minimizing sequence of non-negative functions $(u_k^j)$ which can be
written as $u_k^j=v_k^j+k$ where $v_k^j\in M_0^1(D)$ and $v_k^j$
weakly converges in $M^1_0(D)$ to $v_k$ as $j\to +\infty$ and
$\lim_{j\to +\infty}J(u_k^j)\geq J(u_k)$. Hence $J(u_k)=inf_EJ(u)$.
Hence we have a unique non-negative weak solution $u_k$ of
(\ref{3.2}). By the regularity theory (see \cite{HK} and \cite{JX}),
we can assume that $u_k$ is smooth in $D$, and applying Bony's
maximum principle to $w=-u_k$ (see also \cite{B}) we have $u_k>0$ in
$D$.

We now normalize ${u}_k$ such that $u_k(0)=1$. Then we may use the
Harnack inequality. So, for each compact subset $\Omega$ of $R^n$,
we can find uniform constants $K$ and $C$ such that for each
$k\geq K$, we have
$$
0<u_k(x)\leq C, x\in \Omega.
$$
Using the a priori estimates in \cite{X}, we may assume that there
is a subsequence which  converges in $C^3_{loc}(R^n)$ to a smooth
non-negative solution $u$ of the problem
$$
Hu=\lambda g(x)u,\quad{in} \quad {R^n}.
$$
with $u(0)=1$. Again by Bony's maximum principle, we have $u(x)>0$
for every $x\in R^n$. This proves Theorem \ref{T1.3}. \end{proof}

\section{solutions of a sub-elliptic problem}

Although we shall work only in a bounded regular domain of $R^n$, we
emphasize that in the remaining part of this paper, we always assume
that the vector field group $X_j$, $1\leq j\leq m$, is a family of
Hormander's vector fields on $R^n$.

Assume that $\Omega$ is a bounded regular domain in $R^n$. The main
purpose of this section is to study the existence of a non-trivial
non-negative solution of the following sub-elliptic equation:
$$
   -H u+\lambda a(x)u-b(x)u^p=0 \quad{in}\quad \Omega
$$
with the Dirichlet boundary condition:
$$
u=0\quad{on}\quad \partial\Omega,
$$
where $a(x)$ and $b(x)$ are continuous functions on $\Omega$ with
$\|a\|_\infty<\infty$ and $b(x)$ nonnegative and not identically
zero.
 The basic ingredient in the proof consists of the
following lemma (see \cite{DM}).

\bigskip

\begin{Lem}\label{T4.1} (Comparison Principle) Suppose that $\Omega$ is a bounded
domain in $R^n$, $a(x)$ and $b(x)$ are continuous functions on
$\Omega$ with $\|a\|_\infty<\infty$ and $b(x)$ nonnegative and not
identically zero.
 Let $u_1, u_2\in C^2(\Omega)$
be positive in $\Omega$ and satisfy
\begin{equation}
-Hu_1+ a(x) u_1-b(x)u_1^p\leq 0\leq -Hu_2+a(x)u_2-b(x)u_2^p, \; x\in
\Omega \label{4.1}
\end{equation}
and $\overline{\lim}_{x\to\partial \Omega}(u_2-u_1)\leq 0$, where
$p>1$ is a given constant. Then $u_2\leq u_1$ in $\Omega$.
\end{Lem}

\begin{proof} Let $w_1,w_2$ be $C^2$ nonnegative functions on $\Omega$
vanishing near $\partial \Omega$. Using (\ref{4.1}),  applying
integration by parts and subtracting, we obtain
\begin{align}
&-\int_{\Omega}[(Xu_2)(Xw_2)-(Xu_1)(Xw_1)]\nonumber\\
 &\geq
\int_{\Omega}b (x)[u_2^pw_2-u_1^pw_1]+\int_{\Omega}
a(x)(u_1w_1-u_2w_2). \label{4.2}
\end{align}

For $\epsilon>0$, denote $\epsilon_1=\epsilon$ and
$\epsilon_2=\epsilon/2$ and let
\[v_i=[(u_2+\epsilon_2)^2-(u_1+\epsilon_1)^2]^+/(u_i+\epsilon_i),\; i=1,2.\]
Note $v_i$ can be approximated arbitrarily closely in the $M^1\cap
L^\infty$ norm on $\Omega$ by $C^2$ functions vanishing near
$\partial \Omega$. We see that (\ref{4.2}) holds when $w_i$ is
replaced by $v_i$, $i=1,2$.

Let
\[\Omega_+(\epsilon)=\{x\in\Omega:
u_2(x)+\epsilon_2>u_1(x)+\epsilon_1\}.\] We have that the
integrands
 in (\ref{4.2}) (with $w_i=v_i$) vanish outside
$\Omega_+(\epsilon)$. The integral on the left hand side of
(\ref{4.2}) equals
\[-\int_{\Omega_+(\epsilon)}
|Xu_2-\frac{u_2+\epsilon_2}{u_1+\epsilon_1}Xu_1|^2+
|Xu_1-\frac{u_1+\epsilon_1}{u_2+\epsilon_2}Xu_2|^2,\] which is
non-positive. On the other hand, as $\epsilon\to 0$, the first term
on the right hand side of (\ref{4.2}) converges to
\[\int_{\Omega_+(0)}b(x)(u_2^{p-1}-u_1^{p-1})(u_2^2-u_1^2),\]
while the second term in (\ref{4.2}) converges to 0.  Therefore, we
must have
 $u_1\geq u_2$ on $\Omega$, as required. \end{proof}

Using this lemma, we can prove the main result in this section,
which is the following:

\begin{Prop} \label{T4.2} \ Let $\Omega$ be a bounded regular
domain in $R^n$. Suppose $a(x)$ and $b(x)$ are smooth positive
functions on $\overline{\Omega}$, and let $\mu_1$ denote the first
eigenvalue of $H u=\mu a(x) u$ on $\Omega$ under Dirichlet boundary
condition on $\partial \Omega$. Then the problem
\[Hu=\mu u[a(x)-b(x)u^{p-1}], \quad{in}\quad \Omega\; \]
with the boundary condition:
\[u|_{\partial \Omega}=0\] has a unique positive solution
$u_\mu\in C^1(\Omega)$
for every $\mu>\mu_1$.
\end{Prop}

\begin{proof}

The existence follows from a simple upper and lower solution
argument. Clearly any constant greater than
 or equal to $$M=\max_{\overline{\Omega}}[a(x)/b(x)]^{1/(p-1)}$$
 is an upper solution. Let $\phi$ be a
 positive eigenfunction corresponding to $\mu_1$. Then for each
 fixed $\mu>\mu_1$ and
all small positive $\epsilon$, $\epsilon\phi<M$ is a lower solution.
Thus, by using the monotone method (see Proposition \ref{T5.1}
below) there is at least one positive solution.
 If $u_1$ and $u_2$ are two positive solutions, we apply Lemma
 \ref{T4.1}
 and conclude that $u_1\leq u_2$
and $u_2\leq u_1$ both hold on $\Omega$. Hence $u_1=u_2$. This
proves the uniqueness.
\end{proof}

\section{Poisson Problem and Yamabe Type Equation}

  In this section, we shall study a Yamabe-type or logistic
  equation in the following
  form :
\begin{equation}
Hu+k(x)u-K(x)|u|^{p-1}u=0
\end{equation}
where $k(x)$ and $K(x)$ are given smooth functions in $R^n$, $p>1$
is a constant.

The following result should be well-known to experts and its proof
is similar to that of Theorem 2.10 in \cite{N}.

\begin{Prop}\label{T5.1} Suppose that $U_2\geq U_1$ are entire
super-solution and sub-solution respectively of the equation
\begin{equation}
HU+F(x,U)=0 \label{5.1}
\end{equation}
where $F(x,y)$ is a locally Holder continuous function which is
locally Lipschitz in $y$. Then there is an entire solution $U$ of
(\ref{4.2}) satisfying $U_2\geq U\geq U_1$.
\end{Prop}

We now give the proof of Theorem \ref{T1.4}.

 \begin{proof}
According to our Property (P), we have a
positive constant $\theta_1$ such that\\
 1). The Poisson
equation on $R^n$:
$$
HU(x)=-Cf(x).
$$
has a unique solution $V_{\epsilon}$ with $\lim_{|x|\to
+\infty}V_{\epsilon}(x)=\epsilon$ and $0<V_{\epsilon}(x)\leq
\epsilon$ for any $C\in (0,\theta_1)$ and every $\epsilon\in
(1/3,1/2)$, and \\
 2). The Poisson equation on $R^n$:
$$
HU(x)=Cf(x).
$$
has a unique solution $W_{\epsilon}$ with $\lim_{|x|\to
+\infty}W_{\epsilon}(x)=\epsilon$ and $\epsilon\leq
W_{\epsilon}(x)<1$ for any $C\in (0,\theta_1)$ and every
$\epsilon\in (1/3,1/2)$.

Take $\theta=\theta_1/3$ and $C=2\theta_1/3$. Then, for $K(x)\leq
\theta f(x)$ and $\epsilon\in (1/3,1/2)$, we have
$0<V_{\epsilon}\leq W_{\epsilon}<1$,
\begin{align*}
HV_{\epsilon}+&k(x)V_{\epsilon}-K(x)V_{\epsilon}^p\\
&\leq -Cf(x)+|k(x)|\epsilon+|K(x)|\epsilon^p\leq 0,
\end{align*}
and
\begin{align*}
HW_{\epsilon}+&k(x)W_{\epsilon}-K(x)W_{\epsilon}^p\\
&\geq Cf(x)-|k(x)|\epsilon-|K(x)|\epsilon^p\geq 0.
\end{align*}

By Proposition \ref{T5.1}, we have a solution $U$ of (\ref{5.1})
with $V_{\epsilon}\leq U_{\epsilon}\leq W_{\epsilon}$ and
$\lim_{|x|\to +\infty}U_{\epsilon}(x)=\epsilon$.

\end{proof}

\end{document}